\input amstex
\documentstyle {amsppt}
\pageheight{50.5pc} \pagewidth{32pc}


\topmatter
\author
Guszt\'av Morvai and  Benjamin Weiss 
\endauthor

\title{}
Inferring the Conditional Mean
\endtitle

\rightheadtext{ Inferring the Conditional Mean}

\abstract 
Consider  a stationary  real-valued time series $\{X_n\}_{n=0}^{\infty}$ with a priori   
unknown   distribution. 
The  goal is to   
estimate the conditional expectation  $E(X_{n+1}|X_0,\dots, X_n)$ based on the observations 
$(X_0,\dots, X_n)$ in a pointwise consistent way.   
It is well known that this is not possible at all values of $n$. 
 We will estimate it along  stopping times. 
\endabstract

\address
Guszt\'av Morvai (Corresponding author. Tel.: 36-1-4632867; fax.:36-1-4633147.)
Research Group for Informatics and Electronics 
of the Hungarian Academy of Sciences,   
Budapest, 1521 Goldmann Gy\"orgy t\'er 3, Hungary 
\endaddress

\email morvai\@math.bme.hu
\endemail

\address 
Benjamin Weiss (Tel.: 972-2-658-4388; fax.: 972-2-563-0702.)
Hebrew University of Jerusalem
Jerusalem 91904 Israel
\endaddress

\email 
weiss\@math.huji.ac.il
\endemail

\keywords Nonparametric estimation, stationary processes
\endkeywords

\subjclass 62G05, 60G25, 60G10
\endsubjclass


\endtopmatter

\document
\TagsOnRight


\head Appeared in:  Theory Stoch. Process.  11  (2005),  no. 1-2, 112--120. 
\endhead

 \head Introduction and Statement of Results
\endhead
Suppose  the distribution of the real-valued stationary time series $\{X_n\}_{n=0}^{\infty}$ 
is not known a priori. 
The goal is to estimate 
the conditional expectation $E(X_{n+1}|X_0,\dots,X_n)$ from the data segment 
$X_0,\dots, X_n$ such that 
the difference between the estimate and 
the conditional expectation should tend to zero almost surely as 
the number   of observations  $n$  
tends to infinity. 
This problem (for binary time series) was introduced    
in Cover (1975).   
When  one is obliged to estimate for all $n$,   Bailey (1976)  
and  Ryabko (1988) proved the nonexistence of such a universal  
algorithm  even over  the class of all stationary and ergodic binary time series. 

In a  special case, for certain Gaussian processes,  Sch\"afer  (2002)
constructed an algorithm which can estimate the conditional expectation for every  time 
instance $n$.  

For further reading on related topics cf.  Ornstein (1978), 
Algoet (1992), (1999), 
Morvai Yakowitz and Algoet (1997),  Morvai, Yakowitz and Gy\"orfi 
(1996), Gy\"orfi, Lugosi and Morvai (1999), 
Gy\"orfi and Lugosi (2002), Weiss (2000) and Gy\"orfi et al. (2002).

In this paper we do not require to estimate 
for every time instance $n$, 
but rather, 
merely along  a sequence of stopping times. That is, looking at the data segment 
$X_0,\dots,X_n$ our rule 
will decide if we   estimate for this $n$ or not, 
but anyhow we will definitely estimate 
for infinitely many $n$. 
Algorithms of this kind  were proposed for binary time series in Morvai (2003) and 
Morvai and Weiss (2003).

We will consider 
two-sided real-valued processes $\{X_n\}_{n=-\infty}^{\infty}$. 
A one-sided  stationary time series $\{ X_n\}_{n=0}^{\infty}$ 
can always be considered  to be a two-sided stationary time series 
$\{ X_n\}_{n=-\infty}^{\infty}$.

Let $\Re$ be the set of all real numbers and put ${\Re}^{*-}$  the set of all one-sided 
  sequences of real numbers, 
that is, 
$${\Re}^{*-} =\{ (\dots,x_{-1},x_0): x_i\in \Re \  
\text{for all}  \ -\infty<i\le 0 \}.$$
 Define the metric $d^*(\cdot,\cdot)$ on ${\Re}^{*-}$ as 
$$
d^*((\dots,x_{-1},x_{0}),(\dots,y_{-1},y_{0}))=
\sum_{i=0}^{\infty} 2^{-i-1} {
|x_{-i}-y_{-i}|\over 1+|x_{-i}-y_{-i}| }.
$$ 
 
\bigskip
\noindent
\definition{ Definition:}
The conditional expectation  $E(X_1|\dots,X_{-1},X_{0})$
is almost surely continuous if for some  set $B\subseteq {\Re}^{*-}$  
which has probability one   
the conditional expectation $E(X_1|\dots,X_{-1},X_{0})$ restricted to this  set $B$ 
is continuous with respect to metric $ d^*(\cdot,\cdot)$. 
\enddefinition

\bigskip
\noindent
Now we introduce our algorithm.  
For notational convenience, let 
$X_m^n=(X_m,\dots,X_n)$,
where $m\le n$. Define  the  nested sequence of  partitions 
$\{ {\Cal P}_k\}_{k=0}^{\infty}$ 
of the real line  as follows.
Let 
$$
{\Cal P}_k=\{ [i2^{-k}, (i+1)2^{-k}) \ : \ \text{for}  \ i=0,1,-1,2,-2,\dots \}.
$$   
Let $x\rightarrow [x]^k$ denote a quantizer that assigns to any point $x\in \Re$ 
the  unique interval in ${\Cal P}_k$ that contains $x$. 
Let $[X_m^n]^k=([X_m]^k,\dots,[X_n]^k)$. 

\bigskip
\noindent
We define the stopping times $\{\lambda_n\}$ along which we will estimate. 
Set $\lambda_0=0$. 
For $n=1,2,\ldots$,  define   $\lambda_n$
recursively. 
Let
$$
\lambda_n=\lambda_{n-1}+\min\{t>0 : [X_{t}^{\lambda_{n-1}+t}]^n=
[X_{0}^{\lambda_{n-1}}]^n. \tag1
$$
Note that $\lambda_n\ge n$ and it is a stopping time on $[X_0^{\infty}]^n$. 
Let $f_k:{\Cal P}_k \rightarrow \Re$ denote a function that assigns to any cell 
$A\in {\Cal P}_k$ a point in $A$. 
The $n$th estimate $m_n$ is defined as   
$$
m_n=
{1\over n}\sum_{j=0}^{n-1} f_j([X_{\lambda_j+1}]^j). \tag 2   
$$ 
Observe that $m_n$ depends solely on $[X_{0}^{\lambda_n}]^n$. This estimator can be viewed
as  a sampled version of the predictor in 
Morvai, Yakowitz and Gy\"orfi 
(1996), Weiss (2000), Algoet~(1999) and Gy\"orfi et al. (2002).

\bigskip
\noindent
Define the time series $\{ {\tilde X}_n\}_{n=-\infty}^{0}$ as 
$$
{\tilde X}_{-n}=
\lim_{j\to \infty} X_{\lambda_j-n} \  \text{for}  \ n\ge 0, \tag 3 
$$
 where the limit exists  
 since  the intervals $\{[X_{\lambda_j-n}]^j\}_{j=n}^{\infty}$ 
are nested and their lengths tend to zero.

\bigskip
\noindent 
Define the function $e : {\Re}^{*-}\rightarrow (-\infty,\infty)$ as 
$$e(x^{0}_{-\infty})=
E(X_1|X^{0}_{-\infty}=x^0_{-\infty}).$$
We will prove the following theorem.

\bigskip
\noindent
\proclaim{Theorem}  Let  $\{X_n\}$ be a real-valued stationary time series 
with $E(|X_0|^2)<\infty$. Then almost surely  
$$
\lim_{n\to\infty} m_n=  \lim_{n\to\infty} 
E(X_{\lambda_n+1}|[X_0^{\lambda_n}]^n)= e({\tilde X}^0_{-\infty}) 
$$
and 
$$
\lim_{n\to\infty} \left| m_n- E(X_{\lambda_n+1}|[X_0^{\lambda_n}]^n)\right| =0. 
$$
Moreover, if in addition the conditional expectation  $E(X_1|X_{-\infty}^{0})$ is almost surely 
continuous then almost surely
$$
\lim_{n\to\infty} \left| m_n-E(X_{\lambda_n+1}|X_0^{\lambda_n})\right| =0.
$$  
Unfortunately, there is a stationary and ergodic Markov chain $\{X_n\}$  
 taking values from a countable subset of the unit interval such that
$$
P\left( \limsup_{n\to\infty} \left| m_n-E(X_{\lambda_n+1}|X_0^{\lambda_n})\right| >0\right)>0.
$$ 
\endproclaim

\bigskip
\noindent
\remark{Remarks} 

Let $\{X_n\}$ be a real-valued stationary time series 
with $E(|X_0|^2)<\infty$. 
If the distribution of $X_0$ happens to concetrate on finitely many atoms then 
$$E(X_{\lambda_n+1}|[X_0^{\lambda_n}]^n)=E(X_{\lambda_n+1}|X_0^{\lambda_n})
\ \ \text{eventually}
$$
and so  
$|m_n-E(X_{\lambda_n+1}|X_0^{\lambda_n})|\to 0$ almost surely,
without any continuity condition. 

\bigskip
Let $\{X_n\}$ be a real-valued stationary time series 
with $E(|X_0|^2)<\infty$. 
If one knows in advance that the distribution of $X_0$ concentrates on finite or 
countably infinite atoms 
then one may omit the partition ${\Cal P}_k$, the quantizer $[\cdot]^k$ and the function $f_k(\cdot)$ 
entirely. That is, one may define  $\lambda^{\prime}_0=0$ and  
for $n=1,2,\ldots$ set   
$$
\lambda^{\prime}_n=\lambda^{\prime}_{n-1}+
\min\{t>0 : X_{t}^{\lambda^{\prime}_{n-1}+t}=X_{0}^{\lambda^{\prime}_{n-1}}\}
$$
and     
$$m^{\prime}_n=
{1\over n}\sum_{j=0}^{n-1} X_{\lambda^{\prime}_j+1}.
$$  
Then  
$$\lim_{n\to\infty} \left| m^{\prime}_n-
E(X_{\lambda^{\prime}_n+1}|X_0^{\lambda^{\prime}_n})\right| =0 \ \ \text{almost surely}
$$ 
without any continuity condition. Particularly, $m^{\prime}_n$ works 
for the counterexample process in the third part of the Theorem.      

\bigskip

The counterexample Markov chain in the third part of the Theorem  
of course will not possess almost surely continuous conditional expectation  
$E(X_{1}|X_{-\infty}^0)$. 

\bigskip

From the proof of Bailey (1976), Ryabko (1988), 
Gy\"orfi, Morvai, Yakowitz~(1998) it is clear that    
even for the class of all stationary and ergodic binary time series 
with almost surely continuous conditional expectation 
$E(X_{1}|X_{-\infty}^0)$ 
one can not estimate 
 $E(X_{n+1}|X_0^{n})$ for all $n$ in a pointwise consistent way.  
\endremark

\head 
Proofs
\endhead

It will be useful to define other processes 
 $\{ {\hat X}^{(k)}_n\}_{n=-\infty}^{\infty}$ for $k\ge 0$ as follows.
Let 
$$
\hat X^{(k)}_{-n}=X_{\lambda_k-n} \  \text{for} \  -\infty <n<\infty. \tag 4  
$$

\noindent
For an arbitrary real-valued stationary time series $\{Y_n\}$, 
let $\hat\lambda_0(Y^0_{-\infty})=0$  
and    for $n\ge 1$ define 
$$
\hat\lambda_n(Y^0_{-\infty})=\hat\lambda_{n-1}(Y^0_{-\infty})-
\min\{t>0 : 
[Y_{\hat\lambda_{n-1}-t}^{-t}]^{n}
=
[{Y}_{\hat\lambda_{n-1}}^{0}]^{n}\}.
$$
Let $T$ denote the left shift operator,  
that is, $(T x^{\infty}_{-\infty})_i=x_{i+1}$. It is easy to see that if 
$\lambda_n(x_{-\infty}^{\infty})=l$ then 
${\hat \lambda}_n(T^l x_{-\infty}^{\infty})=-l$.
  
\bigskip
\noindent
\demo{Proof of the Theorem}

\noindent
{\bf Step 1.} 
{\it We show that for arbitrary $k\ge 0$,  
the  time series $\{{\hat X}^{(k)}_n\}_{n=-\infty}^{\infty}$ and  
$\{X_n\}_{n=-\infty}^\infty$ 
have 
identical  distribution.} 

\bigskip
\noindent 
It is enough to show  that  
for all $k\ge 0$, $m\ge n\ge 0$,  and Borel set $F\subseteq {\Re}^{n+1}$,
$$
P(({\hat X}^{(k)}_{m-n},\dots,{\hat X}^{(k)}_{m})\in F)=P(X^m_{m-n}\in F).
$$
 This is immediate by  stationarity of $\{X_n\}$ and by the fact that  
for all $k\ge 0$,  $m\ge n\ge 0$, $l\ge 0$, 
$F\subseteq {\Re}^{n+1}$,
$$
T^{l} \{X^{\lambda_k+m}_{\lambda_k+m-n}\in F,\lambda_k=l\} =
\{ X^{m}_{m-n}\in F,{\hat \lambda}_k(X^0_{-\infty})=-l\}.   
$$

\bigskip
\noindent
{\bf Step 2.} 
{\it We show that for $k\ge 0$, almost surely, 
$$
\hat\lambda_k(\dots,{\hat X}^{(k)}_{-1},{\hat X}^{(k)}_0)=\hat\lambda_k({\tilde X}^0_{-\infty})
$$
and 
$$ 
[ {\tilde X}^0_{\hat\lambda_k({\tilde X}^0_{-\infty})}]^{k+1}=
[{\hat X}^{(k)}_{\hat\lambda_k(\dots,{\hat X}^{(k)}_{-1},{\hat X}^{(k)}_0)},\dots,{\hat X}^{(k)}_0]^{k+1}.
$$
}

\bigskip
\noindent
Since we are dealing with a nested sequence of partitions and 
$\hat\lambda_k$ depends solely on the $k$th quantized sequence,  
 it is enough to prove that for any 
$i\ge 0$ and for all $j\ge i$, almost surely,  
$[{\tilde X}_{-i}]^{j+1}= [\hat X^{(j)}_{-i}]^{j+1}$. 
(Note that $\lambda_j(X_0^{\infty})-j\ge 0$.)
If ${\tilde X}_{-i}\not\in [{\hat X}^{(j)}_{-i}]^{j+1}$ for some $j\ge i$ then  
this must happen at a right  end-point of some  interval in   
$\bigcup_{k=0}^{\infty} {\Cal P}_k$.  
By \thetag{3} and Step~1, we have 
$$
\aligned
 1 &- P( {\tilde X}_{-i}\in [\hat X^{(j)}_{-i}]^{j+1}  \  \text{for all}  \ j\ge i )\\
&\le 
\sum_{k=i}^{\infty}\sum_{s=-\infty}^{\infty}
P({\tilde X_{-i}}=s 2^{-k}, {\hat X}^{(j)}_{-i} <{\tilde X_{-i}} \   \text{for all}
\   j\ge k)\\
&\le \sum_{k=i}^{\infty}\sum_{s=-\infty}^{\infty}
\lim_{j\to\infty} P( s 2^{-k}-2^{-j}\le {\hat X}^{(j)}_{-i} <s 2^{-k})\\
&= \sum_{k=i}^{\infty}\sum_{s=-\infty}^{\infty} \lim_{j\to\infty}
P( s 2^{-k}-2^{-j}\le { X}_{-i} <s 2^{-k})\\
&= 0.
\endaligned
$$

\bigskip
\noindent
{\bf Step 3.}
{\it We show that the distributions of
$\{ {\tilde X}_n\}_{n=-\infty}^{0}$ and $\{ X_n\}_{n=-\infty}^0$ are the same.}

\bigskip
\noindent
This is immediate from  Step~1 and Step~2. 

\bigskip
\noindent
The time series $\{ {\tilde X}_n\}_{n=-\infty}^{0}$ 
is stationary, since  $\{ X_n\}_{n=-\infty}^0$ is stationary, and it can be 
extended to be a two-sided time series
$\{ {\tilde X}_n\}_{n=-\infty}^{\infty}$. 
We will use this fact only for the purpose of defining the conditional expectation 
$E({\tilde X}_1|{\tilde X}^{0}_{-\infty})$.   

\bigskip
\noindent
{\bf Step 4.}
{\it We prove the first part of the Theorem.} 

\bigskip
\noindent
Consider 
$$
\align
m_n &=
{1\over n}\sum_{j=0}^{n-1} 
\left(  
f_j([X_{\lambda_j+1}]^j)  -E(f_j([X_{\lambda_j+1}]^j)|[X_{0}^{\lambda_j}]^j) \right)\\ 
&+  
{1\over n}\sum_{j=0}^{n-1} \left(E(f_j([X_{\lambda_j+1}]^j)|[X_{0}^{\lambda_j}]^j) -
E(X_{\lambda_j+1}|[X_{0}^{\lambda_j}]^j)\right)\\
&+
{1\over n}\sum_{j=0}^{n-1} E(X_{\lambda_j+1}|[X_{0}^{\lambda_j}]^j). \tag 5
\endalign
$$

\noindent
Observe that
 $\{\Gamma_j=f_j([X_{\lambda_j+1}]^j)  -E(f_j([X_{\lambda_j+1}]^j)|[X_{0}^{\lambda_j}]^j) \}$ 
is a sequence of orthogonal random variables with $E\Gamma_j=0$ and     
$E\left( \Gamma_j^2\right) \le E\left( |X_1|^2\right)+2E|X_1|+1$ since  
$E\left( \Gamma_j^2\right) \le E\left( |X_{\lambda_j+1}|^2\right)+2E|X_{\lambda_j+1}|+1$
and, by Step~1, $X_{\lambda_j+1}$ has the same distribution as $X_1$. 
Now by Theorem~3.2.2 in R\'ev\'esz (1968),
$$
{1\over n}\sum_{j=0}^{n-1} \Gamma_j \to 0 \  \text{almost surely.}
$$

\noindent
The second term tends to zero since $|f_j([X_{\lambda_j+1}]^j)-X_{\lambda_j+1}|\le2^{-j}$.
\noindent
Now we deal with the third term.
By Step~2, Step~1 and Step~3, 
$$
E(X_{\lambda_j+1}|[X_{0}^{\lambda_j}]^j)=
E({\tilde X}_{1}|[{\tilde X}_{\hat\lambda_j({\tilde X}^0_{-\infty})}^0]^j).
$$ 
The latter forms a martingale and by Theorem 7.6.2 in Ash~(1972), almost surely,   
$$
E(X_{\lambda_j+1}|[X_{0}^{\lambda_j}]^j)=
E({\tilde X}_{1}|[{\tilde X}_{\hat\lambda_j({\tilde X}^0_{-\infty})}^0]^j)\to 
E({\tilde X}_{1}|{\tilde X}_{-\infty}^0). \tag 6
$$
By ~\thetag{5} and ~\thetag{6}, almost surely, 
$$
\lim_{n\to\infty} m_n= E({\tilde X}_{1}|{\tilde X}_{-\infty}^0). \tag 7
$$
Thus the first part of the Theorem is proved. 

\bigskip
\noindent
{\bf Step 5.}
{\it We prove the second part of the Theorem. }

\bigskip
\noindent
By~\thetag{7} it is enough to prove that
almost surely  
$E(X_{\lambda_j+1}|X_{0}^{\lambda_j})\to E({\tilde X}_{1}|{\tilde X}_{-\infty}^0)$ provided that 
$E(X_1| X_{-\infty}^0)$ is almost surely continuous. 
By assumption, the function $e(\cdot)$ is continuous on a set $B\subseteq {\Cal \Re}^{*-}$ 
with $P(X^0_{-\infty}\in B)=1$. 
By~Step~1
  and Step~3,  
$$
P({\tilde X}^{0}_{-\infty}\in B, (\dots,{\hat X}^{(j)}_{-1},{\hat X}^{(j)}_0)\in B 
\  \text{for all}  \  j\ge0)=1. \tag 8
$$

\noindent 
Let
$$
{\Cal  N}_j(X_0^{\lambda_j})=\{ z^0_{-\infty}\in \Re^{*-} \ : \ 
z_{-\lambda_j}\in [X_{0}]^j,\dots, 
z_{0}\in [X_{\lambda_j}]^j \}.   
$$
By  \thetag{4}, \thetag{8} and Step~2, almost surely, for all $j$,  
$$
(\dots,{\hat X}_{-1}^{(j)}, {\hat X}_{0}^{(j)})\in {\Cal  N}_j(X_0^{\lambda_j})\bigcap B \  
\text{and} \  
{\tilde X}^{0}_{-\infty} \in {\Cal  N}_j(X_0^{\lambda_j})\bigcap B. \tag 9
$$
Put
$$
\Theta_j(X_0^{\lambda_j})=\sup_{y^0_{-\infty}, z^0_{-\infty}\in {\Cal  N}_j(X_0^{\lambda_j})
\bigcap B}
|e(y^0_{-\infty})-e(z^0_{-\infty})|.
$$
Since $e(\cdot)$ is continuous  on set $B$ and 
by \thetag{9}, almost surely, 
$$
\lim_{j\to\infty} 	\Theta_j(X_0^{\lambda_j})=0. \tag 10
$$
By~ \thetag{9} and~ \thetag{10}, almost surely, 
$$
\align
\limsup_{j\to\infty}&\left|E\left( e({\tilde X}^{0}_{-\infty})|[X_0^{\lambda_j}]^j\right)- 
E\left( e(\dots,{\hat X}_{-1}^{(j)}, {\hat X}_{0}^{(j)})
|X_0^{\lambda_j}\right)\right|\\
&\le
\limsup_{j\to\infty} 
E \left(\left| E\left( e({\tilde X}^{0}_{-\infty} )|[X_0^{\lambda_j}]^j\right)
- 
e(\dots,{\hat X}_{-1}^{(j)}, {\hat X}_{0}^{(j)}) \right| |X_0^{\lambda_j}\right)\\
&\le 
 \limsup_{j\to\infty} 
E\left( 	\Theta_j(X_0^{\lambda_j})|X_0^{\lambda_j}\right)\\
&=    \limsup_{j\to\infty} 	\Theta_j(X_0^{\lambda_j})\\
&= 0.\tag 11
\endalign
$$
By~Step~2, 
$$
\align
E\left(X_{\lambda_j+1}|X_0^{\lambda_j}\right)&=
E\left( e({\tilde X}^{0}_{-\infty})|[{\tilde X}^0_{{\hat \lambda}_j}]^j \right)\\
&-
\left\{ E\left( e({\tilde X}^{0}_{-\infty})|[X_0^{\lambda_j}]^j \right)-
E\left( e(\dots,{\hat X}_{-1}^{(j)}, {\hat X}_{0}^{(j)})
|X_0^{\lambda_j}\right)\right\}.
\endalign
$$
The first term tends to $e({\tilde X}^{0}_{-\infty})$ by the 
almost sure martingale convergence theorem (cf.Theorem 7.6.2 in Ash~(1972))
since by ~Step~3,  
$E \left| e({\tilde X}^{0}_{-\infty})\right| \le E\left| {\tilde X}_1\right| =
 E\left| X_1\right| <\infty$. 
The second term tends to zero by ~\thetag{11}. 
The proof of the second part of the  Theorem is complete.

\bigskip
\noindent
{\bf Step 6.}
{\it We prove the third part of the Theorem.} 

\bigskip
\noindent 
First we define a Markov chain $\{M_n\}$ on the nonnegative integers which will 
serve as a technical tool for our counterexample 
process. 
Let the transition probabilities be as follows. 
$$P(M_1=0|M_0=0)=P(M_1=1|M_0=0)=P(M_1=0|M_0=1)=2^{-1}$$ 
and for $i=2,3,\dots, $ let 
$$P(M_1=i|M_0=1)=2^{-i} \  \text{and}  \  P(M_1=0|M_0=i)=1.$$ 
All other transitions happen with probability zero. Note that one can reach state $1$ only from state $0$. 
It is easy to see that the Markov chain just defined yields a  stationary and ergodic 
time series with 
initial probabilities $P(M_0=0)={4\over 7}$, $P(M_0=1)={2\over 7}$, and for $i=2,3,\dots$ 
$P(M_0=i)={1\over 7} {1 \over  2^{i-1} }$. 
Our counterexample process $\{X_n\}$ will be a one to one  function of 
the Markov chain $\{M_n\}$. 
Define the function $h: \{0,1,2,\dots\} \rightarrow \Re$  as 
$h(0)=0$, $h(1)=1$ and for $i\ge 2$ put $h(i)={2^{-2^i}\over 2}$. 
Let $X_n=h(M_n)$. 
Since $h(\cdot)$ is one to one, $\{X_n\}$ is also a Markov chain. Since $\{{\tilde X}_n\}$ 
has the same distribution as $\{X_n\}$, 
$\{{\tilde X}_n\}$ is also a Markov chain. 
Let 
$$
A_n=\{ h(i): h(i)<2^{-(n+1)} \ \text{for}  \ i=0,1,2,\dots \}.
$$
Note that $h(i)\in A_n$ if and only if $[h(i)]^{n+1}=[0]^{n+1}$. 
Define the event 
$$
H=\{ {\tilde X}_0=0, X_0^1=(0,1)\}.
$$
Observe: If $X_1=1$ then $X_0=0$. (State $1$ can be reached only from state $0$.)
The event $\{{\tilde X}_0=0\}$ happens if and only if 
$ X_{\lambda_n} \in A_n$ for all $n=1,2,\dots$. Since $[h(0)]^1=[h(i)]^1$ for $i\ge 2$ and 
for all $k\ge 0$,  $[h(1)]^k\neq [h(i)]^k$ provided $i\neq 1$
the event 
$\{{\tilde X}_{-1}=1\} $ occurs if and only if $X_1=1$. 
It follows that 
$$
H=\{ X_0=0, X_1=1, X_{\lambda_n} \in A_n \  \text{for}  \   n=1,2,\dots \}=
\{ {\tilde X}_{-2}^0=(0,1,0)  \}.
$$
Since the time series $\{ {\tilde X}_n\}$ has the same distribution as $\{X_n\}$, 
$$
P(H)=P( X_{-2}^0=(0, 1, 0) )={4 \over 7} {1\over 2} {1\over 2}={1\over 7}>0.
$$
It will be enough to show that $X_{\lambda_n}\in  A_n-\{0\}$ happens infinitely often given the condition $H$ since 
if $X_{\lambda_n}\in  A_n-\{0\}$ happens then $X_{\lambda_n+1}=0$ and 
by ~\thetag{7}, on $H$ 
$$m_n\to E({\tilde X}_{1}| {\tilde X}_{0}=0)=0.5 $$
and so 
 $$
P\left( \limsup_{n\to \infty} |m_n- E(X_{\lambda_n+1}|X_0^{\lambda_n})|=0.5| H\right)=1
 $$ and $P(H)>0$. To prove that $\{X_{\lambda_n}\in  A_n-\{0\}\}$ occurs infinitely often we need the following observation for repeated use: 
By the Markov property and the construction in~\thetag{1} 
if  $x_i\in A_i$ for $i=1,2,\dots, j$ then for $j\ge 1$, 
$$
P(X_{\lambda_j}=x_j|X_0^1=(0,1), X_{\lambda_m}=x_m \  \text{for}  \ 1\le m<j) = 
P(X_1=x_j|X_0=1,X_1\in A_{j-1}). \tag 12 
$$
Indeed, for $j=1$ this is trivial, since $X_1=1$ implies that $X_0=0$, $\lambda_1=2$ 
while $X_0=1$ implies that $X_1\in A_0$.  For 
$j\ge 2$ set $\psi_0^j=\lambda_{j-1}-1$ and for $i\ge 1$ the $\psi_i^j$ will be the successive occurrences of the block
$[X_0^{\lambda_{j-1}-1}]^{j}$
in the $j$-th quantization, defined by  
$$
\psi_i^j=
\min\{t>\psi_{i-1}^j : [X_{t-\lambda_{j-1}+1}^{t}]^{j}=
[X_{\psi_{i-1}^j-\lambda_{j-1}+1}^{\psi_{i-1}^j}]^{j}\}. 
$$
These $\psi_i^j$ are   stopping times for $i=1,2,\dots$.
Temporarily let $D_j$ denote the event 
$$\{X_0^1=(0,1), X_{\lambda_m}=x_m \  \text{for} \  1\le m<j\}.$$ 
The way that $\lambda_j$ is defined means that on $D_j$ if $\lambda_j$ occurs at the $i$-th repetition of 
$[X_0^{\lambda_{j-1}-1}]^{j}$ it is because $\psi_i^j<\lambda_j$ and $X_{\psi_i^j+1} \in A_{j-1}$. It follows that
$$
P(X_{\lambda_j}=x_j|D_j) =
\sum_{i=1}^{\infty}  
P(X_{\psi_i^j+1}=x_j|X_{\psi_i^j+1}\in A_{j-1}, \psi_i^j<\lambda_j,D_j) P( \psi_i^j+1=\lambda_j | D_j). 
$$
Since $x_j\in A_j\subseteq A_{j-1}$, each expression 
$P(X_{\psi_i^j+1}=x_j|X_{\psi_i^j+1}\in A_{j-1}, \psi_i^j<\lambda_j,D_j)$ can be written  as 
$$
P(X_{\psi_i^j+1}=x_j|X_{\psi_i^j+1}\in A_{j-1}, \psi_i^j<\lambda_j,D_j)={P(X_{\psi_i^j+1}=x_j| \psi_i^j<\lambda_j,D_j)
\over 
P(X_{\psi_i^j+1}\in A_{j-1}| \psi_i^j<\lambda_j,D_j)}
$$
and then by decomposition according to the value $l$ of $\psi_i^j$ we get 
$$
\align
&P(X_{\psi_i^j+1}=x_j| \psi_i^j<\lambda_j,D_j) \\
&=
\sum_{l=1}^{\infty}  \left(
{
P(X_{l+1}=x_j|\psi_i^j=l<\lambda_j, D_j)\over 
P(X_{l+1}\in A_{j-1} |\psi_i^j=l<\lambda_j, D_j) } 
P(\psi_i^j=l, X_{\psi_i^j+1}\in A_{j-1}|\psi_i^j<\lambda_j,  D_j)\right). 
\endalign
$$
Observe that $X_{\psi_i^j}=1$ provided $X_1=1$ and the event $\{\psi_i^j<\lambda_j\}$ is 
measurable with respect to 
$\sigma([X_0^{\psi_i^j}]^j)$.  
 Now by the Markov property we get 
 $$
 \align
 &P(X_{\psi_i^j+1}=x_j|X_{\psi_i^j+1}\in A_{j-1}, \psi_i^j<\lambda_j,D_j)\\
 &=
 \sum_{l=1}^{\infty} \left( 
{
P(X_{l+1}=x_j|X_l=1)\over 
P(X_{l+1}\in A_{j-1} |X_l=1) } \cdot 
{
P(\psi_i^j=l, X_{\psi_i^j+1}\in A_{j-1}|\psi_i^j<\lambda_j, D_j)\over 
 P( X_{\psi_i^j+1}\in A_{j-1}|\psi_i^j<\lambda_j, D_j)} \right).
 \endalign
 $$
 By stationarity and since $x_j\in A_j\subseteq A_{j-1}$,   
$$
 {P(X_{l+1}=x_j|X_l=1)\over 
P(X_{l+1}\in A_{j-1} |X_l=1) } =P(X_1=x_j| X_{1}\in A_{j-1}, X_0=1) .
$$
Combining all this we get
$$
\align
&P(X_{\lambda_j}=x_j| D_j) \\
&=
 P(X_1=x_j| X_{1}\in A_{j-1}, X_0=1)  \\
 &\cdot \left(\sum_{i=1}^{\infty} 
P( \psi_i^j+1=\lambda_j |  D_j) 
\sum_{l=1}^{\infty} 
{
P(\psi_i^j=l, X_{\psi_i^j+1}\in A_{j-1}| \psi_i^j<\lambda_j, D_j)\over 
 P( X_{\psi_i^j+1}\in A_{j-1}| \psi_i^j<\lambda_j, D_j)} \right)  \\
&= 
P(X_1=x_j| X_{1}\in A_{j-1}, X_0=1)  \sum_{i=1}^{\infty}  
P( \psi_i^j+1=\lambda_j |  D_j)  \\
&= P(X_1=x_j| X_{1}\in A_{j-1}, X_0=1)
\endalign
$$
and we have proved \thetag{12}.

\noindent
In order to show that  the events $$\{X_{\lambda_n}\in  A_n-\{0\} \}$$ occur 
infinitely often we prove that 
 they 
have sufficiently large  conditional probabilities and they are conditionally independent given the condition $H$. 
First we calculate $P(X_{\lambda_n}\in  A_n-\{0\}|H)$. For $n\ge 2$, by~\thetag{12},   
$$
\align
&P(X_{\lambda_n}\in  A_n-\{0\}|H)\\
&= { P(\{X_{\lambda_n}\in  A_n-\{0\}\} \bigcap H)\over P(H)  }\\
&=    
{ P(X_{\lambda_n}\in  A_n-\{0\}|X_0^1=(0,1), X_{\lambda_j}\in A_j \  
\text{for} \  1\le j<n) \over 
P(X_{\lambda_n}\in A_n |X_0^1=(0,1), X_{\lambda_j}\in A_j \  \text{for}  \ 1\le j<n) }\\
&\cdot
\prod_{m=n+1}^{\infty} {
P(X_{\lambda_m}\in  A_m|X_0^1=(0,1),X_{\lambda_n}\in  A_n-\{0\},  
X_{\lambda_j}\in A_j \  \text{for} \  1\le j<m) \over 
 P(X_{\lambda_m}\in  A_m, |X_0^1=(0,1), X_{\lambda_j}\in A_j \ 
\text{for} \  1\le j<m)   } \\
&=  { P(X_{\lambda_n}\in  A_n-\{0\}|X_0^1=(0,1), X_{\lambda_j}\in A_j \ 
\text{for} \  1\le j<n) \over 
P(X_{\lambda_n}\in A_n |X_0^1=(0,1), X_{\lambda_j}\in A_j \  \text{for}  \ 1\le j<n) }\\
&\ge  P(X_{\lambda_n}\in  A_n-\{0\}|X_0^1=(0,1), X_{\lambda_j}\in A_j \ 
\text{for} \  1\le j<n) \\
&= P(X_1\in A_n, X_1\neq 0|X_0=1, X_1\in A_{n-1} )\\
&\ge P(X_1\in A_n, X_1\neq 0|X_0=1)\\
&=  \sum_{i\in A_n-\{0\} } {1\over 2^{i} }\\
&= \sum_{i>\log_2 (n) } {1\over 2^i}\\
&\ge  {1\over n}. 
\endalign
$$
We have just proved that 
$$ 
\sum_{n} P(X_{\lambda_n}\in  A_n-\{0\}|H)\ge \sum_n {1\over n} =\infty. \tag{13}
$$
Now we will prove that for $n=1,2,\dots$, the events $\{X_{\lambda_n}\in  A_n-\{0\} \}$ are conditionally independent given $H$.
Since 
$$
\align
&P(X_{\lambda_i}\in A_i-\{0\} \  \text{for}  \ i=1,2,\dots,k|H)\\
&=
\sum_{x_1\in A_i-\{0\} }\dots \sum_{x_k\in A_k-\{0\} } P(X_{\lambda_i}=x_i \  \text{for}  \ i=1,2,\dots,k|H)
\endalign
$$
it is enough to show that the events $\{X_{\lambda_i}=x_i\}$ are conditionally independent given the condition $H$,
 provided that $x_i\in A_i$. 
Let $x_i\in A_i$. Then by repeated use of ~\thetag{12}
$$
\align
&P(X_{\lambda_i}=x_i \  \text{for}  \ i=1,2,\dots,k|H)\\
&= {P(X_{\lambda_i}=x_i \  \text{for}   \ i=1,2,\dots,k , H) \over P(H)} \\
&=  
\left( \prod_{m=1}^k 
{ P(X_{\lambda_m}=x_m| X_0^1=(0,1), X_{\lambda_j}= x_j \ \text{for}  \  1\le j<m) \over 
P(X_{\lambda_m}\in A_m| X_0^1=(0,1),
 X_{\lambda_j}\in A_j \  \text{for}  \ 1\le j<m)}\right)   \\
&\cdot \prod_{l=k+1}^{\infty}  { P(X_{\lambda_l}\in A_l|X_0^1=(0,1), X_{\lambda_i}=x_i  \ 
\text{for} \ 1\le i\le k \  \text{and}  \  X_{\lambda_j}\in A_j \  \text{for}  \
 1\le j<l) \over 
P(X_{\lambda_l}\in A_l| X_0^1=(0,1),  
 X_{\lambda_j}\in A_j   \  \text{for}  \ 1\le j<l) } \\
&=
 \prod_{m=1}^k 
{ P(X_{\lambda_m}=x_m| X_0^1=(0,1), X_{\lambda_j}\in A_j  \   \text{for}  \ 1\le j<m) \over 
P(X_{\lambda_m}\in A_m| X_0^1=(0,1),
 X_{\lambda_j}\in A_j \  \text{for}  \ 1\le j<m)}   \\
&= 
\prod_{m=1}^k \left( 
{ P(X_{\lambda_m}=x_m| X_0^1=(0,1), 
X_{\lambda_j}\in A_j  \ \text{for}  \ 1\le j<m) \over 
P(X_{\lambda_m}\in A_m| X_0^1=(0,1),
 X_{\lambda_j}\in A_j \  \text{for}  \ 1\le j<m)} \right.\\
&\cdot    \left.
\prod_{l=m+1}^{\infty}  { P(X_{\lambda_l}\in A_l|X_0^1=(0,1), 
X_{\lambda_j}\in A_j  \   \text{for}  \  1\le j<l) \over 
P(X_{\lambda_l}\in A_l|X_0^1=(0,1), X_{\lambda_j}\in A_j  \   \text{for}  \
1\le j<l) }\right) \\
&=
\prod_{m=1}^k \left( 
{ P(X_{\lambda_m}=x_m| X_0^1=(0,1), X_{\lambda_j}\in A_j \   \text{for}  \ 
1\le j<m) \over 
P(X_{\lambda_m}\in A_m| X_0^1=(0,1),
 X_{\lambda_j}\in A_j \  \text{for}  \ 1\le j<m)} \right.\\
&\cdot    \left.
\prod_{l=m+1}^{\infty}  { P(X_{\lambda_l}\in A_l|X_0^1=(0,1),X_{\lambda_m}=x_m,  
X_{\lambda_j}\in A_j  \  \text{for}  \ 1\le j<l) \over 
P(X_{\lambda_l}\in A_l|X_0^1=(0,1), 
X_{\lambda_j}\in A_j \    \text{for}  \ 1\le j<l) }\right) \\
&= \prod_{i=1}^k {P(X_{\lambda_i}=x_i,H)\over P(H) }\\
&= \prod_{i=1}^k P(X_{\lambda_i}=x_i|H).
\endalign
$$
Now by \thetag{13} and the Borel-Cantelli lemma 
(cf. Lemma B in R\'enyi (1970) on page 390) the events $\{X_{\lambda_n}\in  A_n-\{0\} \}$ occur infinitely often and 
the third part of the Theorem is proved.  
The proof of the Theorem is complete.

\enddemo

\enddocument